\documentclass{amsart}
\usepackage{graphicx}
\usepackage{amscd}
\usepackage{amsmath}
\usepackage{amsfonts}
\usepackage{amssymb}
\newtheorem{theorem}{Theorem}

\newtheorem{remark}[theorem]{Remark}

\numberwithin{equation}{section}

\begin{document}
\title[Isoperimetry and Symmetrization for Sobolev spaces]{Isoperimetry and Symmetrization for Sobolev spaces on metric spaces}
\author{Joaquim Martin$^{\ast}$}
\address{Department of Mathematics\\
Universitat Aut\`onoma de Barcelona}
\email{jmartin@mat.uab.cat}
\author{Mario Milman}
\address{Department of Mathematics\\
Florida Atlantic University}
\email{extrapol@bellsouth.net}
\urladdr{http://www.math.fau.edu/milman}
\thanks{2000 Mathematics Subject Classification Primary: 46E30, 26D10.}
\thanks{$^{\ast}$ Supported by Grants MTM2007-60500 and by 2005SGR00556.}

\begin{abstract} Using isoperimetry we obtain
new symmetrization inequalities that allow us to provide a unified framework
to study Sobolev inequalities in metric spaces. The applications include
concentration inequalities, as well as metric versions of the P\'{o}%
lya-Szeg\"{o} and
Faber-Krahn principles.
\end{abstract}
\maketitle
\section{Introduction}

\label{}

This is a follow up to our recent work \cite{mami0}, where we obtained new
symmetrization inequalities for Sobolev functions that compare the
rearrangement of a function with the rearrangement of its gradient, and
incorporate in their formulation the isoperimetric profile (cf. (\ref{rea})
below). These inequalities imply in a straightforward fashion functional
inequalities for very general rearrangement invariant norms or quasi-norms
(e.g. $L^{p},$ Orlicz, Lorentz, Marcinkiewicz spaces). One remarkable
characteristic of these inequalities is that they preserve their form as we
move from one measure space to another, the only thing that changes are the
corresponding isoperimetric profiles. As a consequence we were able to provide
a unified framework to study the classical Sobolev-Poincar\'{e} inequalities,
logarithmic Sobolev inequalities, as well as concentration inequalities (cf.
\cite{LE} and the references therein). Importantly, if the isoperimetric
profile does not depend on the dimension (like in the Gaussian case) then the
corresponding inequalities are dimension free.

The purpose of this note is to outline the modifications that are necessary to
extend our earlier results to the setting of metric spaces. Indeed, under
relatively weak assumptions, all the tools that we need are available in the
metric setting (cf. \cite{BH}), and our methods can be readily adapted to
provide an almost painless extension. In particular, the results of this note,
when combined with the method developed\footnote{Our method builds on a
variant of Maz'ya's truncation principle, combined with the relevant
isoperimetric inequalities, the co-area formula and classical arguments from
real interpolation theory (cf. Calder\'{o}n \cite{CA}). We call this method to
obtain symmetrization inequalities ``symmetrization via truncation''.} in
\cite{mami0}, produce concentration inequalities in metric spaces, as well as
a Sobolev metric space version of the P\'{o}lya-Szeg\"{o} principle; while our
results combined with the method of [\cite{KM}, Theorem 3] imply metric
Faber-Krahn inequalities.

Let $\left(  \Omega,d,\mu\right)  $ be a metric space equipped with a
separable Borel probability measure $\mu$. For measurable functions
$u:\Omega\rightarrow\mathbb{R},$ the distribution function of $u$ is given by
$\lambda_{u}(t)=\mu(\{x\in{\Omega}:\left|  f(x)\right|  >t\}$ $(t>0)$, the
decreasing rearrangement $u^{\ast}$ of $u$ is defined, as usual, by $u^{\ast
}(s)=\inf\{t\geq0:\lambda_{u}(t)\leq s\}$ $(t\in(0,1)]),$ and we let
$u^{\ast\ast}(t)=\frac{1}{t}\int_{0}^{t}u^{\ast}(s)ds.$ For $A\subset\Omega,$
a Borel set, let $Per(A)=\lim\inf_{\varepsilon\rightarrow0}\frac{\mu\left(
A_{\varepsilon,d}\right)  -\mu\left(  A\right)  }{\varepsilon},$ where
$A_{\varepsilon,d}=\left\{  x\in\Omega:\exists y\in A\mbox{
}d(x,y)<\varepsilon\right\}  $ denotes the $\varepsilon-$extension of $A$ with
respect to the metric $d.$ An isoperimetric inequality measures the relation
between $Per(A)$ and $\mu(A)$ by means of the isoperimetric profile
$I=I_{(\Omega;d;\mu)}$, which is defined as the pointwise maximal function
$I:[0,1]\rightarrow\left[  0,\infty\right)  $, such that $Per(A)\geq
I(\mu(A)),$for all Borel sets $A$. Finally, in this setting for a given
Lipschitz function $f$ (we shall write in what follows $f\in Lip(\Omega)$) the
modulus of the gradient is defined, as usual, by $|\nabla f(x)|=\lim
\sup_{d(x,y)\rightarrow0}\frac{|f(x)-f(y)|}{d(x,y)}.$

\section{Main results}

\begin{theorem}
\label{teomain}Suppose that the isoperimetric profile $I$ is concave,
continuous, increasing on $(0,1/2)$ and symmetric about the point $1/2.$ Then
the following statements hold\footnote{except where indicated all
rearrangements are with respect to the measure $\mu$.} (and in fact are
equivalent):
\[
(i):\text{ }\forall\text{ }f\in Lip(\Omega),\text{ }\int_{0}^{\infty}%
I(\lambda_{f}(s))ds\leq\int_{{\Omega}}\left|  \nabla f(x)\right|
d\mu(x)\text{ (Ledoux)}.
\]%
\[
(ii):\text{ }\forall\text{ }f\in Lip(\Omega),\text{ }(-f^{\ast})^{\prime
}(s)I(s)\leq\frac{d}{ds}\int_{\{\left|  f\right|  >f^{\ast}(s)\}}\left|
\nabla f(x)\right|  d\mu(x)\text{ (Talenti-Maz'ya).}%
\]%
\[
(iii):\text{ }\forall\text{ }f\in Lip(\Omega),\text{ }\int_{0}^{t}((-f^{\ast
})^{\prime}(.)I(.))^{\ast}(s)ds\leq\int_{0}^{t}\left|  \nabla f\right|
^{\ast}(s)ds\text{ (P\'{o}lya-Szeg\"{o}).}%
\]
(The second rearrangement on the left hand side is with respect to the
Lebesgue measure).
\begin{equation}
(iv):\text{ }\forall\text{ }f\in Lip(\Omega),\text{ }(f^{\ast\ast}(t)-f^{\ast
}(t))\leq\frac{t}{I(t)}\left|  \nabla f\right|  ^{\ast\ast}(t).\label{rea}%
\end{equation}
\end{theorem}

Given any rearrangement invariant space\footnote{A Banach lattice of functions
$X(\Omega)$ is called a rearrangement invariant (r.i.) space (cf. \cite{BS})
if $g\in X(\Omega)$ implies that all functions $f$ \ with the same decreasing
rearrangement, $f^{\ast}=g^{\ast},$ also belong to $X(\Omega),$ and, moreover,
$\Vert f\Vert_{X(\Omega)}=\Vert g\Vert_{X(\Omega)}$. There is an essentially
unique r.i. space $\bar{X}(0,1)$ of functions on the interval $(0,1)$
consisting of all $g:(0,1)\rightarrow R$ such that $g^{\ast}(t)=f^{\ast}(t)$
for some function $f\in X(\Omega)$.} $X(\Omega)$, it follows readily from
(\ref{rea}) that for all Lip functions, we have
\begin{equation}
\left\|  f\right\|  _{LS(X)}:=\left\|  \left(  f^{\ast\ast}(t)-f^{\ast
}(t)\right)  \frac{I(t)}{t}\right\|  _{\bar{X}}\leq\left\|  \nabla f\right\|
_{X}.\label{corres}%
\end{equation}

\begin{remark}
For $L^{1}$ norms these Poincar\'{e} inequalities are a simple variant of
Ledoux's inequality $(i)$. Indeed, let $m_{f}$ be a median\footnote{i.e.
$\mu\left(  f\geq m\right)  \geq1/2$ and $\mu\left(  f\leq m\right)  \geq
1/2$.} of $f$, then it is easy to see that
\begin{equation}
\int_{{\Omega}}\left|  f-m_{f}\right|  d\mu\leq\frac{1}{2I(1/2)}\int_{{\Omega
}}\left|  \nabla f(x)\right|  d\mu(x).\label{pr1}%
\end{equation}
\end{remark}

The novelty of our inequalities, and the corresponding associated spaces
$LS(X),$ is that they incorporate the isoperimetric profiles associated with
the geometry in question. These spaces are not necessarily normed, although
often they are equivalent to normed spaces (cf. \cite{PU}), and, in the
classical cases, lead to optimal Sobolev-Poincar\'{e} inequalities (cf.
\cite{MMP}, \cite{mami0}, \cite{MM1} and the references therein).

We now investigate the optimality of the Poincar\'{e} type inequality
(\ref{corres}). The following result is new in the context of r.i. spaces.

\begin{theorem}
\label{optimal}Let $(\Omega,\mu)=(R^{n},\mu_{r}^{\otimes n}),$ with $\mu
_{r}(x)=$ $\varphi(x)dx,I_{\mu_{r}^{\otimes n}}(t)\approx\varphi
(F^{-1}(t)),\;t\in\left[  0,1\right]  $, where $F^{-1}$ is the inverse of the
distribution function associated to the density $\varphi(x)dx$\footnote{This
choice of $I$ is motivated by the results in \cite{Bo}, \cite{Bob} and
\cite{BFR}.}. Let $X(\Omega)$, $Y(\Omega)$ be r.i. spaces. Then, the following
statements are equivalent
\begin{equation}
(i):\text{ }\forall\text{ }f\in Lip(\Omega),\,\left\|  f-\int f\right\|
_{Y}\preceq\left\|  \nabla f\right\|  _{X}.\label{poin}%
\end{equation}%
\[
(ii):\text{ }\left\|  \int_{t}^{1}f(s)\frac{ds}{I(s)}\right\|  _{\bar{Y}%
}\preceq\left\|  f\right\|  _{\bar{X}},\text{ }\forall\text{ }0\leq f\in
\bar{X},\text{with }supp(f)\subset(0,1/2).
\]
Moreover,

\begin{enumerate}
\item [(a)]If the operator $Q_{I\text{ }}f(t)=\frac{I(t)}{t}\int_{t}%
^{1}f(s)\frac{ds}{I(s)}$ is bounded from $\bar{X}$ to $\bar{X},$ then the next
inequality can be added to the list of equivalences
\begin{equation}
(iii):\text{ }\left\|  f\right\|  _{\bar{Y}}\preceq\left\|  f^{\ast}%
(t)\frac{I(t)}{t}\right\|  _{\bar{X}}.\label{perdida}%
\end{equation}

\item[(b)] On the other hand if $Q_{I\text{ }}$ is not bounded from $\bar{X} $
to $\bar{X},\ $but $\left\|  f\right\|  _{X}\simeq\left\|  f^{\ast\ast
}\right\|  _{\bar{X}}$, then the next inequality can be added to the list of
equivalences
\begin{equation}
\left\|  f\right\|  _{\bar{Y}}\preceq\left\|  f\right\|  _{LS(X)}+\left\|
f\right\|  _{L^{1}.}\label{harhar}%
\end{equation}
\end{enumerate}
\end{theorem}

As a concrete illustration\footnote{For further examples we refer to
\cite{BH}, \cite{MiE}, and the references therein.} consider the family of
probability measures on the real line given by $d\mu_{r}(t)=\alpha_{r}%
^{-1}e^{-\left|  t\right|  ^{r}}dt=\varphi_{r}(t)dt$, $1<r\leq2,$ where
$\alpha_{r}^{-1}$ is chosen to ensure that $\mu_{r}(\mathbb{R)=}1.$ These
probabilities form a \ scale between exponential and Gaussian measure. The
associated isoperimetric profile is given by $I_{\mu_{r}}(t)=\varphi_{r}%
(F_{r}^{-1}(t)),$ where $F_{r}^{-1}$ is the inverse of the distribution
function associated to the density $\varphi_{r}(t)$ (cf. \cite{Bov})$.$ The
isoperimetric profiles $I_{\mu_{r}^{\otimes n}},$ associated to the product
probability measures $\mu_{r}^{\otimes n},$ is dimension free (see
\cite{BFR}): there is a universal constant $c(r)$ such that $I_{\mu_{r}%
}(t)\geq\inf_{n\geq1}I_{\mu_{r}^{\otimes n}}(t)\geq c(r)I_{\mu_{r}}(t)$. As an
application let $n\geq2$, and apply Theorem \ref{optimal} to $X=L^{p}%
(\mathbb{R}^{n},d\mu_{r}^{\otimes n}),$ $1\leq p<\infty,$ then (cf. also
\cite[Theorem 3]{MM1}),
\[
\int_{0}^{1}\left(  \left(  f-\int f\right)  ^{\ast}(s)\frac{I_{\mu_{r}}%
(s)}{s}\right)  ^{p}ds\preceq\int_{\mathbb{R}^{n}}\left|  \nabla f(x)\right|
^{p}d\mu_{r}^{\otimes n}(x),
\]
with dimension free constants. In particular, since (see \cite[Lemma
16.1]{Bov}) $\lim_{t\rightarrow0^{+}}\frac{I_{\mu_{r}}(t)}{t\left(  \log
\frac{1}{t}\right)  ^{1/q}}=r$ , $1/r+1/q=1$, it follows easily that
\[
\int_{0}^{1}f^{\ast}(s)^{p}(\log\frac{1}{s})^{p/q}ds\preceq\int_{\mathbb{R}%
^{n}}\left|  \nabla f(x)\right|  ^{p}d\mu_{r}^{\otimes n}(x)+\int
_{\mathbb{R}^{n}}\left|  f(x)\right|  ^{p}d\mu_{r}^{\otimes n}(x).
\]
Moreover, for this class of measures, $L^{p}(LogL)^{p/q}$ is the best possible
choice among all r.i. spaces $Y$ for which the inequality $\Vert f-\int
f\Vert_{Y}\preceq\Vert\nabla f\Vert_{L^{p}}$ holds. If $p=\infty,$ we have
\begin{equation}
\Vert f-\int f\Vert_{LS(L^{\infty})}=\left\|  \left(  \left(  f-\int f\right)
^{\ast\ast}(t)-\left(  f-\int f\right)  ^{\ast}(t)\right)  \frac{I_{\mu_{r}%
}(t)}{t}\right\|  _{L^{\infty}}\preceq\left\|  \nabla f\right\|  _{L^{\infty}%
}.\label{concen}%
\end{equation}
The relation to concentration inequalities follows directly from our main
inequality. Indeed, we have
\[
\sup_{t<1}\left\{  (f^{\ast\ast}(t)-f^{\ast}(t))\frac{I_{\mu_{r}}(t)}%
{t}\right\}  \leq\sup_{t}\left|  \nabla f\right|  ^{\ast\ast}(t)=\left\|
f\right\|  _{Lip},
\]
which, by the asymptotic properties of $I_{\mu_{r}},$ implies that
$f^{\ast\ast}(t)-f^{\ast}(t)\preceq\frac{\left\|  f\right\|  _{Lip}}{\left(
\log\frac{1}{t}\right)  ^{1/q}}$ ($0<t<1/2$). We may now proceed as in
[\cite{mami0}, Section 7].

Let us finally consider Sobolev embeddings into $L^{\infty}.$ Notice that from
inequality (\ref{rea}) we get
\[
\left\|  f\right\|  _{\infty}-2\int_{0}^{1/2}f^{\ast}(t)=\int_{0}%
^{1/2}(f^{\ast\ast}(t)-f^{\ast}(t))\frac{dt}{t}\leq\int_{0}^{1/2}\left(
\frac{1}{t}\int_{0}^{t}\left|  \nabla f\right|  ^{\ast}(s)ds\right)  \frac
{dt}{I_{\mu_{r}}(t)}=\int_{0}^{1/2}\left|  \nabla f\right|  ^{\ast}(s)\int
_{s}^{1/2}\frac{ds}{I_{\mu_{r}}(s)s}.
\]
Using the asymptotics of $I_{\mu_{r}}(s)$ combined with the Poincar\'{e}
inequality \ref{pr1} yields
\[
\left\|  f-m_{f}\right\|  _{\infty}\preceq\int_{0}^{1/2}\left|  \nabla
f\right|  ^{\ast}(s)\frac{ds}{s\left(  \log\frac{1}{s}\right)  ^{1/q}}.
\]

\section*{Acknowledgements}

We are very grateful to the referee for suggestions to improve the
presentation.The first named author has been partially supported by Grants
MTM2007-60500, MTM2008-05561-C02-02 and by 2005SGR00556.
%etc, etc
%
%
%
%
%
%
%

%The Appendices part is started with the command \appendix;
% \appendix
%
%
%
%
%
%
%

%\section{}
% \label{}
%
%
%
%
%
%
%

%The Acknowledgements are an un-numbered section
%\section*{Acknowledgements}
% Acknowledgements text here
%
%
%
%
%
%
%

\end{document}